\documentclass[preprint,12pt]{elsarticle}

\usepackage{amssymb}
\usepackage{amsmath}
\usepackage{algorithm}
\usepackage{algpseudocode}

\usepackage{float} % for [H] option
\newcommand{\mc}{MC}
\newcommand{\agac}{AG\&AC}
\newcommand{\mcac}{MC\&AC}
\newcommand{\ircm}{IRCM}
\newcommand{\mcircm}{MC\&IRCM}

\usepackage{hyperref}
\usepackage{cleveref}
\newcommand{\kitaichi}[2]{\bar{\chi}({#1},{#2})}
\newcommand{\kitaichisum}[2]{{\bar{\chi}}_{\scalebox{0.5}{\text{MC}}} ({#1},{#2})}

\begin{document}

\begin{frontmatter}



\title{Chromatic number of random graphs: an approach using a recurrence relation}


% Author name         
\author[1]{Yayoi Abe\corref{cor1}}
\ead{yayoi.abe@akane.waseda.jp}
\cortext[cor1]{Corresponding author}

\author{Auna Setoh}
\ead{keruberosu@asagi.waseda.jp}

\author[1]{Gen Yoneda}
\ead{yoneda@waseda.jp}

\affiliation[1]{organization={Graduate School of Fundamental Science and Engineering, Waseda University},
             addressline={3-4-1 Okubo},
             city={Shinjuku},
             postcode={169-8555},
             state={Tokyo},
             country={Japan}}

%% Abstract
\begin{abstract}
The vertex coloring problem to find chromatic numbers is known to be unsolvable in polynomial time. Although various algorithms have been proposed to efficiently compute chromatic numbers, they tend to take an enormous amount of time for large graphs. In this paper, we propose a recurrence relation to rapidly obtain the expected value of the chromatic number of random graphs. Then we compare the results obtained using this recurrence relation with other methods using an exact investigation of all graphs, the Monte Carlo method, the iterated random color matching  method, and the method presented in Bollob\'{a}s' previous studies. 
\end{abstract}



%% Keywords

\end{frontmatter}



\section{Introduction}
The vertex coloring problem is one of the graph coloring problems. The problem is to find the minimum number of colors used in a coloring assignment such that no two adjacent vertices of a graph have the same color. Such coloring is called proper (vertex) coloring, and the minimum number of colors required for proper coloring is called the chromatic number. Finding the chromatic number is known to be %\ac{nph}%
nondeterministic polynomial-time hard (NP-hard)\cite{garey1979computers} and solving the vertex coloring problem is very difficult. Therefore, many heuristic algorithms have been developed\cite{10.1093/comjnl/10.1.85, laguna2001grasp, hertz1987using}. Such algorithms enable the efficient computation of chromatic numbers without investigating all coloring patterns. However, it still tends to take an enormous amount of time to determine the chromatic number of a very large graph.
  
Meanwhile, a different approach, the probabilistic method, has been used to examine the properties of random graphs.  Bollob\'{a}s\cite{alma991002309759704032} proved that the chromatic number $\chi(G_{n})$ satisfies   
\begin{equation}
\frac{n}{2 \log_d n} \left( 1+\frac{\log_d \log n}{\log n} \right)\le \chi(G_{n}) \le \frac{n}{2 \log_d n} \left( 1+\frac{3\log_d \log n}{\log n} \right), 
\label{eq:bollobas1}
\end{equation}
where $d=1/(1-p)$, in which $p$ is the probability that an edge exists between two vertices, and $G_n$ is a random graph with $n$ vertices.

Bollob\'{a}s also showed that the chromatic number of almost every random graph is equal to the lower bound of \eqref{eq:bollobas1} as $n \to \infty$\cite{bollobas1988chromatic}. However, because \eqref{eq:bollobas1} is an estimation for a very large graph, his estimate is not applicable to small graphs. Furthermore, the range between the upper and lower bounds of \eqref{eq:bollobas1} is too wide for the estimation of the expected value of the chromatic number.
  
In this paper, we propose a stochastic recurrence relation to find the expected value of the chromatic number of binomial random graphs with $n$ vertices and the probability $p$, where $p$ is the probability that an edge exists between two vertices. The total number of patterns of simple labeled graphs with $n$ vertices is $2^{\binom{n}{2}}$. We calculate the expected value of the chromatic number of these graphs using the stochastic recurrence relation and compare the expected value with the results from the heuristic algorithm and Bollob\'{a}s' bounds of \eqref{eq:bollobas1}.

%In Sect.~\ref{sec:Prelim}, 
In \Cref{sec:Prelim}, 
we give the definitions, notations, and exactly obtained chromatic numbers. \Cref{sec:MC} contains a description of the experimentally obtained expected value of the chromatic number. In \Cref{sec:Rec}, we propose our stochastic recurrence relation that predicts the expected value of the chromatic number of random graphs. Finally, we conclude in \Cref{sec:Con}. Some results of computational experiments are reported in Appendix.
%\Cref{sec:A,sec:B,sec:C}.

\section{Preliminaries}
\label{sec:Prelim}
We provide the definitions and notations for random graphs in \Cref{subsec:Pre-def} and describe how to obtain the exact chromatic numbers in \Cref{subsec:Pre-exact}. 

\subsection{Definitions and notations}
\label{subsec:Pre-def}

  In general, a probability space $(\Omega,P)$ consists of a sample space $\Omega$ and a probability function $P : \Omega \to [0,1]$ satisfying $\sum_{\omega \in \Omega} P(\omega)=1$. 
We consider a probability space whose sample space is a set consisting of all $2^{\binom{n}{2}}$ simple labeled graphs with the vertex set $V=[n]=\{ 1,2,\dots,n \}$. This sample space is equivalent to a set consisting of all spanning subgraphs of the complete graph $K_n$. This set is denoted as $\mathcal{G}_n$. Each element $G_n$ of the sample space $\mathcal{G}_n$ is called a random graph. The number of vertices in a graph is called \textit{order} and is denoted by the letter $n$. An \textit{n-order graph} is a graph with $n$ vertices.

Given the probability $p \in [0,1]$, each possible edge in the random graph $G_n$ is connected independently to the probability $p$. Then, the probability function $P_p (G_{n}), G_{n} \in \mathcal{G}_{n}$ is 
\begin{equation}
P_p (G_{n}) = p^{|E(G_{n})|} (1-p)^{\binom{n}{2}-|E(G_{n})|}, 
\end{equation}
where $|E(G_{n})|$ is the number of edges in the graph $G_{n}$ and $0 \le |E(G_{n})| \le \binom{n}{2}$. Throughout this paper, $p$ is fixed and independent of $n$. The probability function $P_p(G_n)$ satisfies $\sum_{G_n \in \mathcal{G}_n} P_p(G_n)=1$ \cite{bollobas1998modern,Frieze_Karoński_2023}. 
This probability space is denoted as $(\mathcal{G}_n, P_p)$. The random graph $G_{n} \in \mathcal{G}_{n}$, which is selected according to the probability function $P_p(G_n)$, is called a \textit{binomial random graph}. This model was introduced by Gilbert\cite{gilbert1959}.

Our aim is to calculate the expectation of the chromatic number of  a random graph $G_n \in \mathcal{G}_n$. We denote this expectation as $\kitaichi{n}{p}$. Let $\chi(G_n)$ be the chromatic number of the random graph $G_n \in \mathcal{G}_n$, then $\kitaichi{n}{p}$ is written as 
\begin{equation}
\kitaichi{n}{p} =\sum_{G_{n} \in \mathcal{G}_n} \chi(G_{n}) P_p (G_{n}) 
= \sum_{G_{n} \in \mathcal{G}_n} \chi(G_{n}) \, p^{|E(G_{n})|} (1-p)^{\binom{n}{2}-|E(G_{n})|}.
\label{eq:chi-bar}
\end{equation}  
In particular, if $p=1/2$, then
\begin{equation}
\kitaichi{n}{1/2}=\sum_{G_{n} \in \mathcal{G}_n} \frac{1}{2^{\binom{n}{2}}} \, \chi(G_{n}). 
\end{equation}

\begin{figure}[htbp]
 \centering
\includegraphics[width=0.8\linewidth]{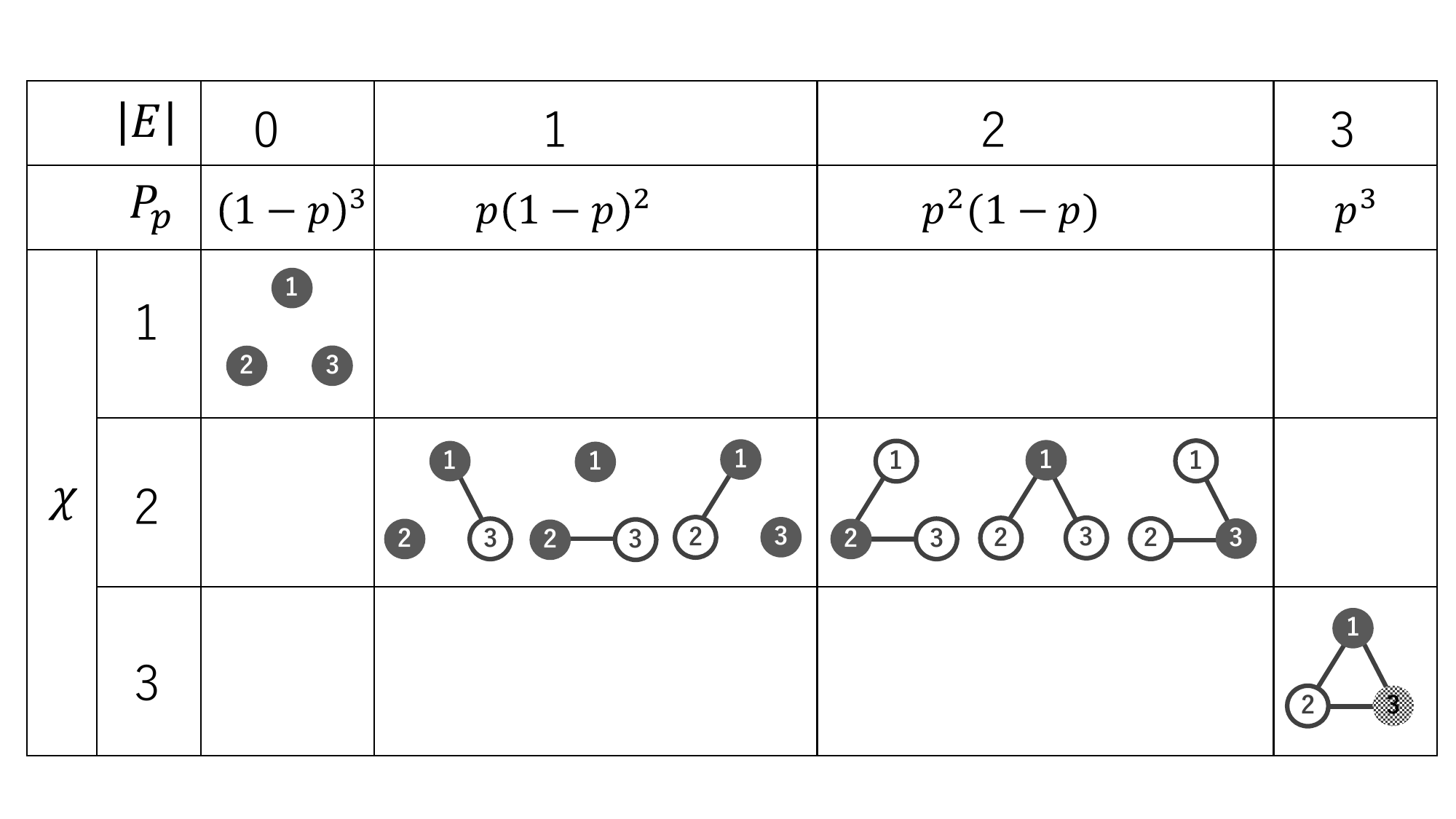} 
 \caption{Probability space $(\mathcal{G}_3, P_p)$. The top line  is the number of edges $|E(G_3)|$ and the second line is the probability function $P_p(G_3)$. }
 \label{fig:fig_G3p}
\end{figure}

For example, the elements of the sample space in $\mathcal{G}_{3}$ are the eight spanning subgraphs shown in \Cref{fig:fig_G3p}. The chromatic number $\chi(G_{3})$ and the probability function $P_p(G_{3})$ are also displayed in \Cref{fig:fig_G3p}. 
The expected value of $\chi(G_{3})$, $\kitaichi{3}{p}$, is
\begin{equation}
\begin{split}
\kitaichi{3}{p} &= 1 \times (1-p)^3 + 2 \times 3p(1-p)^2 
+ 2 \times 3p^2 (1-p) + 3 \times p^3 \\
&= 2 p^3-3 p^2+3 p+1. 
\label{eq:bar-chi-3p}
\end{split}
\end{equation}
\subsection{Exact chromatic numbers of random graphs}
\label{subsec:Pre-exact}
For graphs with a small order, all coloring patterns can be examined for all graph patterns. We call this survey method
%\ac{agac}.
all graphs and all coloring search (\agac). 
In \agac, for each graph of order $n$, we first examine whether each coloring pattern in $n^n$ possible coloring patterns is proper coloring. Next, we seek for the minimum number of colors used in proper coloring.  Then, we get the minimum number as the chromatic number of the graph. We perform this for all $2^{\binom{n}{2}}$ graph patterns. 
Therefore, the expected value obtained by \agac \ corresponds to the exact value.   
  
We were able to examine the graphs up to order $n=9$ using \agac \ in realistic time\footnote{In our study, we define ``realistic time'' as ``within 20 days''. } on our computer. 
The results obtained by \agac \ for graphs with order $n=2$ to $9$ are shown in \Cref{tb:exact_edge2,tb:exact_edge3,tb:exact_edge4,tb:exact_edge5,tb:exact_edge6,tb:exact_edge7,tb:exact_edge8,tb:exact_edge9} in \Cref{sec:A}. These tables show the number of graphs with the edge number $|E| \,(0 \le |E| \le 2^{\binom{n}{2}})$ and the chromatic number $\chi \,(1 \le \chi \le n)$. 

For $n=1$, the element of $\mathcal{G}_1$ is only one graph with no edges. Therefore, $\kitaichi{1}{p}$ is one. 
For $n=2$ to $9$, the expected value of the chromatic number is expressed in polynomial form in the same way that polynomial \eqref{eq:bar-chi-3p} was obtained. The expected values $\kitaichi{n}{p}$, $n=2,\dots,9$ are expressed as polynomials of the probability $p$ as follows.   

\begin{align}
\kitaichi{1}{p} =& 1 \\
\kitaichi{2}{p} =& 1+p \\
\kitaichi{3}{p} =& 2 p^3-3 p^2+3 p+1 \\
\kitaichi{4}{p} =& 3 p^6-15 p^4+24 p^3-15 p^2+6 p+1 
\end{align}
\begin{equation}
\begin{split}
\kitaichi{5}{p} =& -10 p^{10}+90 p^9-240 p^8+310 p^7-265 p^6 \\
   &+234 p^5-210 p^4 +130 p^3-45 p^2+10 p+1
\end{split}
\end{equation}
\begin{equation}
\begin{split}
\kitaichi{6}{p} =& -175 p^{15}+1470 p^{14}
     -5475 p^{13}+12240 p^{12} \\
     &-17775 p^{11} +15960 p^{10}-6230 p^9 -3645 p^8 \\
     &+7035 p^7 -5405 p^6+2985 p^5-1365 p^4 \\
     &+475 p^3-105 p^2+15 p+1 
\end{split}
\end{equation}
\begin{equation}
\begin{split}
\kitaichi{7}{p} =& -949 p^{21}+4725 p^{20}
    +14175 p^{19}-202195 p^{18} \\
    &+880215 p^{17}-2275014 p^{16}
    +3997770 p^{15} -5036475 p^{14} \\
    &+4648770 p^{13}-3157945 p^{12}
    +1586382 p^{11}-639786 p^{10} \\
    &+290885 p^9-187215 p^8
    +116850 p^7-55769 p^6 \\
    &+20391 p^5-5985 p^4+1365 p^3
    -210 p^2+21 p+1
\end{split}
\end{equation}
\begin{equation}
\begin{split}
\kitaichi{8}{p} =& 29701 p^{28}-615412 p^{27}
     +5615022 p^{26}-30989560 p^{25}
     +118474685 p^{24} \\
     &-339285492 p^{23} +767083590 p^{22}
     -1420223760 p^{21}+2203036185 p^{20} \\
     &-2889645360 p^{19}+3195414040 p^{18} 
     -2942158968 p^{17}+2210870585 p^{16} \\
     &-1315890324 p^{15}+588213840 p^{14}
     -172019372 p^{13} +11727835 p^{12} \\
     &+20100528 p^{11}-14697998 p^{10}
     +7145320 p^9-3067617 p^8 +1182440 p^7 \\
     &-380940 p^6
     +98532 p^5-20475 p^4+3332 p^3-378 p^2+28 p+1
\end{split}
\end{equation}
\begin{equation}
\begin{split}
\kitaichi{9}{p} =& 756370 p^{36}-17953236 p^{35}
      +198489942 p^{34}-1362517968 p^{33} \\
      &+6486125310 p^{32}-22481416692 p^{31}
      +57249093300 p^{30} \\
      &-102158239248 p^{29} 
      +97781480685 p^{28}+88988260340 p^{27} \\
      &-621340032978 p^{26} 
      +1562700838860 p^{25} -2744668858953 p^{24} \\
      &+3777632482680 p^{23} 
      -4251933252108 p^{22} +3990846997140 p^{21} \\
      &-3149579225331 p^{20} +2091295480932 p^{19}-1161021072498  p^{18} \\
      &+531668425428 p^{17} 
      -196512115905 p^{16}+57042006564 p^{15} \\
      &-12991012668 p^{14}+2874704364 p^{13} 
      -1079354325 p^{12} \\
      &+565290936 p^{11} 
      -257330430 p^{10}+95465300 p^9 \\
      &-30222918 p^8+8339220 p^7
      -1957578 p^6+377748 p^5 \\
      &-58905 p^4+7224 p^3-630 p^2+36 p+1 
\end{split}
\end{equation}

\section{Estimation of chromatic number using Monte Carlo method}
\label{sec:MC}

 In the previous section, we dealt with the exact chromatic number of small graphs. 
In this section, we focus on the situation where the graph is so large that an exact chromatic number cannot be obtained. For such graphs, the number of patterns becomes so large that the time required to investigate them is enormous. To reduce computing time, we use the Monte Carlo (\mc) method and experimentally predict the expected value of the chromatic number. 
 
The samples are generated using the \mc \ method as follows. First, starting with an empty graph with $n$ vertices, we select two vertices and connect them with an edge with the probability $p$ (for example, if $p=0.3$, we generate a random number from 0 to 9 and connect the vertices with an edge only when the number is $0, 1$, or $2$). By performing this independently for each $\binom{n}{2}$ vertex pair,  we can generate a graph. This process is repeated until the required number of samples is obtained.

The chromatic number of each graph pattern is obtained either by examining all possible coloring patterns discussed in \Cref{subsec:MC-AC} or by using the iterated random color matching method discussed in \Cref{subsec:MC-IRCM}. Note that from this point on, the discussion is about experimentally predicted values, not exact values.

\subsection{Monte Carlo method and all coloring search}
\label{subsec:MC-AC}
For graphs of relatively small order of even more than $9$, all possible coloring patterns can be searched within a realistic time frame. The chromatic number of each graph obtained by examining all coloring patterns is an exact value. However, using the \mc \ method, the expected value is a prediction, because it is obtained by examining not all graph patterns, but some of the extracted graph patterns. This method is denoted as the Monte Carlo method and all coloring search (\mcac). 

  Let $s$ be the number of samples extracted from the probability space $(\mathcal{G}_n, P_p)$ by the \mc \ method, and let $\{ \chi_1,\chi_2,\dots,\chi_{s} \}$ be the chromatic number data for each sample, then the sample average is%
\begin{equation}
  \kitaichisum{n}{p} = \frac{1}{s} \sum_{i=1}^{s} \chi_i .
\end{equation}
From the law of large numbers, $\lim_{s \to \infty} \kitaichisum{n}{p} = \kitaichi{n}{p}$. The results from \mcac \ for the sample number $s=16384$ are shown in \Cref{tb:monte_all} in \Cref{sec:B}, which are the expected values we were able to obtain in realistic time. 

\subsection{Monte Carlo method and iterated random color matching method}
%\subsection{\ac{mcircm}}
\label{subsec:MC-IRCM}

Although the \mc \ method reduces the number of random graphs searched, the number of possible coloring patterns increases significantly in large graphs; thus, an enormous amount of time is required to find the exact chromatic number for each random graph. To reduce the computation time, we use an algorithm that efficiently finds the chromatic number as close to the exact value as possible.

The algorithm used in this study is the iterated random color matching (\ircm) method. This algorithm is a type of local search method\footnote{In general, a local search method is a heuristic method to find an exact solution or an approximate solution close to an exact solution by iterating operations to improve the obtained solution. }, and is shown as \Cref{algo:R2C}. 
In the \ircm \ algorithm, the chromatic number is obtained by the following procedure: First, every vertex is in a different color. Next, two vertices are selected randomly. Then, one of them is changed to the same color as the other vertex. If any adjacent vertex of the changed vertex has the same color, it is restored to its original color. If the color of one vertex is not changed, do the same with the other vertex. This ``random color matching'' is iterated sufficiently many times.

The number of colors used for coloring may decrease but never increase with each iteration. In most cases, the number of colors is expected to asymptotically approach the exact chromatic number after a sufficient number of iterations. However, this is not true for all cases. In some cases, depending on the sequence in which the vertex pairs are selected, the number of colors may remain larger than the exact chromatic number even after additional iterations.

The method using both the \mc \ method and this IRCM method is denoted as \mcircm. In this study, iteration is terminated when the number of colors in all extracted graphs does not change even after doubling the number of iterations. \Cref{tb:monte_IRCM} in \Cref{sec:C} shows the results obtained by \mcircm \ for 16384 samples with probabilities $p=0.3, 0.5$, and $0.7$.
  
\begin{algorithm}[tb]
\caption{Iterated random color matching method (\ircm)}
\label{algo:R2C}
\begin{algorithmic}[1]
    \renewcommand{\algorithmicrequire}{\textbf{Input:}}
    \renewcommand{\algorithmicensure}{\textbf{Output:}}
    \Require {A graph $G$ of order $n$, maximum number of iterations $t_{max}$}
    \State Create an initial coloring by choosing a different color for each vertex $v_i \in V$.
    \State $t \gets 1$ 
    \While {$t \le t_{max}$} 
        \State Choose $v_i,v_j \in V, i \ne j$, uniformly at random.
        \State Let $c(v)$ be the color number of $v$. 
        \If {$c(v_i) \ne c(v_j)$}
        	\If {the color number of all vertices adjacent to $v_i$ is not $c(v_j)$}
          			\State $c(v_i) \gets c(v_j)$
        	\ElsIf {the color number of all vertices adjacent to $v_j$ is not $c(v_i)$}
          		\State $c(v_j) \gets c(v_i)$
        	\EndIf
      \EndIf
      \State $t \gets t+1$ 
    \EndWhile
    \Ensure Number of colors used in graph $G$
\end{algorithmic}
\end{algorithm}

\subsection{Results}
\label{sec:exact_result} 
\begin{figure}[tb]
 \centering
 \includegraphics[width=0.9\linewidth]{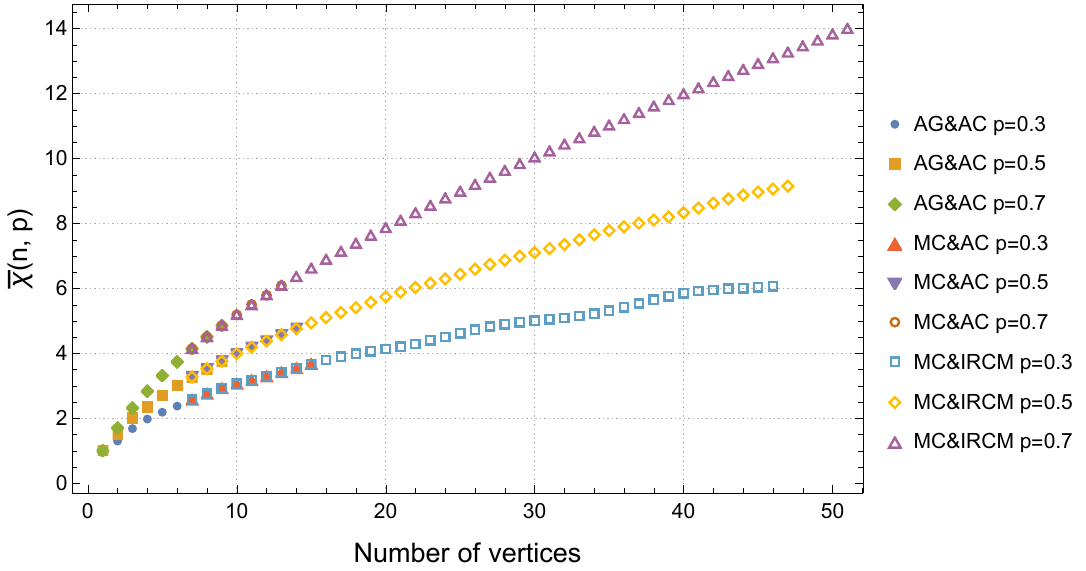} 
 \caption[comparison with expected values of chromatic number]{Expected values of the chromatic number $\kitaichi{n}{p}$ obtained by \agac, \mcac \ and \mcircm. The number of samples for the \mc \ method is 16384. }
 \label{fig:fig_exact}
\end{figure}

\Cref{fig:fig_exact} presents the expected values of the  chromatic number obtained by \agac, \mcac, and \mcircm \ for $p=0.3, 0.5$, and $0.7$, where the number of samples in the \mc \ method was 16384.

By comparing the values obtained by the \mc \ method with the exact values for $n=7, 8$, and $9$, we found that the values obtained by \mcac \ and \mcircm \ were almost identical to the exact values, which were obtained by \agac. Considering this, even for $n>9$, the experimental values obtained by \mcac \ and \mcircm \ are probably very close to the exact values. Thus, we compared the expected values obtained using our new method with these experimental values.

\section{Stochastic recurrence relation}
\label{sec:Rec}
For large graphs with $n>50$, the time required to examine their chromatic number would even more enormous, even if the \mcircm \ method were used. To predict the expected value in a shorter time, we propose a stochastic recurrence relation. We first explain the method of estimation, followed by the results of using the recurrence relation. 
  
\subsection{Estimation method}
 A stochastic recurrence relation is derived as follows. 
Suppose a graph of order $n$ and its expected value of chromatic number $\kitaichi{n}{p}$ are given. A new $(n+1)$-order graph is created by adding one vertex to the graph and connecting this new vertex to $k$ vertices, which are selected from the original graph, where $0\le k \le n$. Let the probability of $\kitaichi{n+1}{p}=\kitaichi{n}{p}$ be $q_n$. Then, the probability of $\kitaichi{n+1}{p}=\kitaichi{n}{p}+1$ is $1-q_n$. Then, 
\begin{equation}
\kitaichi{n+1}{p}=q_n \, \kitaichi{n}{p} + (1-q_n) (\kitaichi{n}{p} + 1) = \kitaichi{n}{p} + 1 - q_n.
\label{eq:exp_rec1}
\end{equation}

The original $n$-order graph is assumed to be colored with $\kitaichi{n}{p}$ colors on average. Then, the created $(n+1)$-order graph is colored with $\kitaichi{n+1}{p}=\kitaichi{n}{p}+1$ colors if all $\kitaichi{n}{p}$ colors are already used at vertices adjacent to the newly added vertex, whereas it is colored with $\kitaichi{n+1}{p}=\kitaichi{n}{p}$ colors if $\kitaichi{n}{p}-1$ or fewer colors are used at vertices adjacent to the new vertex.

 Next, $q_n$ is estimated as follows. The number of vertices connected to the newly added vertex is assumed to be $np$ on average. If $\kitaichi{n+1}{p}=\kitaichi{n}{p}$, the number of colors used at $np$ vertices is at most  $\kitaichi{n}{p}-1$ colors. Then, the probability that one of the $\kitaichi{n}{p}$ colors is not used at all $np$ vertices is $((\kitaichi{n}{p}-1)/\kitaichi{n}{p})^{np}$, and the number of choices for a color that is not used is $\kitaichi{n}{p}$, $q_n$ is approximated as 
\begin{equation}
q_n=\kitaichi{n}{p} \left( \frac{\kitaichi{n}{p} -1}{\kitaichi{n}{p}} \right)^{np}.
\label{eq:exp_rec2}
\end{equation}
Therefore, a recurrence relation about the expected value of the chromatic number is from \eqref{eq:exp_rec1} and  \eqref{eq:exp_rec2}, 
\begin{equation}
\kitaichi{n+1}{p}= \kitaichi{n}{p} + 1 - \kitaichi{n}{p} \left( \frac{\kitaichi{n}{p} -1}{\kitaichi{n}{p}} \right)^{np}.
\label{eq:exp_rec3}
\end{equation}
Although we cannot solve this recurrence relation analytically, we can solve it numerically. The maximum order of the graph for which the chromatic number could be obtained using \mcircm \ was limited to about 50, but using the recurrence relation, we can obtain the expected values even for 100,000 vertices. If the estimation using \eqref{eq:exp_rec3} is accurate, our proposed recurrence relation is very useful.

\subsection{Result and discussion}
In this section, the expected value of the chromatic number  obtained using the recurrence relation \eqref{eq:exp_rec3} is compared with the results presented in \Cref{sec:exact_result} and with the results of the studies by Bollob\'{a}s. As shown in \eqref{eq:bollobas1}, for an $n$-order random graph $G_{n} \in \mathcal{G}_{n}$ with the probability $p$, Bollob\'{a}s proved that the chromatic number $\chi(G_n)$ satisfies 
\begin{equation*}
\frac{n}{2 \log_d n} \left( 1+\frac{\log_d \log n}{\log n} \right)\le \chi(G_{n}) \le \frac{n}{2 \log_d n} \left( 1+\frac{3\log_d \log n}{\log n} \right),
\end{equation*}
where $d=1/(1-p)$. Furthermore, Bollob\'{a}s proved that almost every random graph approaches the lower bound of \eqref{eq:bollobas1} when $n \to \infty$\cite{alma991002309759704032}. 

For the value of the initial term of the recurrence relation \eqref{eq:exp_rec3}, two cases are considered: the case of the exact value $\kitaichi{9}{p}$ and the case of the predicted value on graphs of the largest order examined by \mcircm.
Because we assume that the number of vertices is sufficiently large in the recurrence relation \eqref{eq:exp_rec3}, the initial term should be the value at the largest possible number of vertices. The largest order for which exact values were obtained was $9$. Thus, if an exact value is used for the first term, $\kitaichi{9}{p}$ should be used. Similarly, if the value obtained by \mcircm \ is used for the first term, the value with the largest order among those obtained by \mcircm \ should be used. 

\subsubsection{Comparison of results between \mcircm \ and the recurrence relation methods}

\begin{figure}[H]
 \centering
 \includegraphics[width=0.8\linewidth]{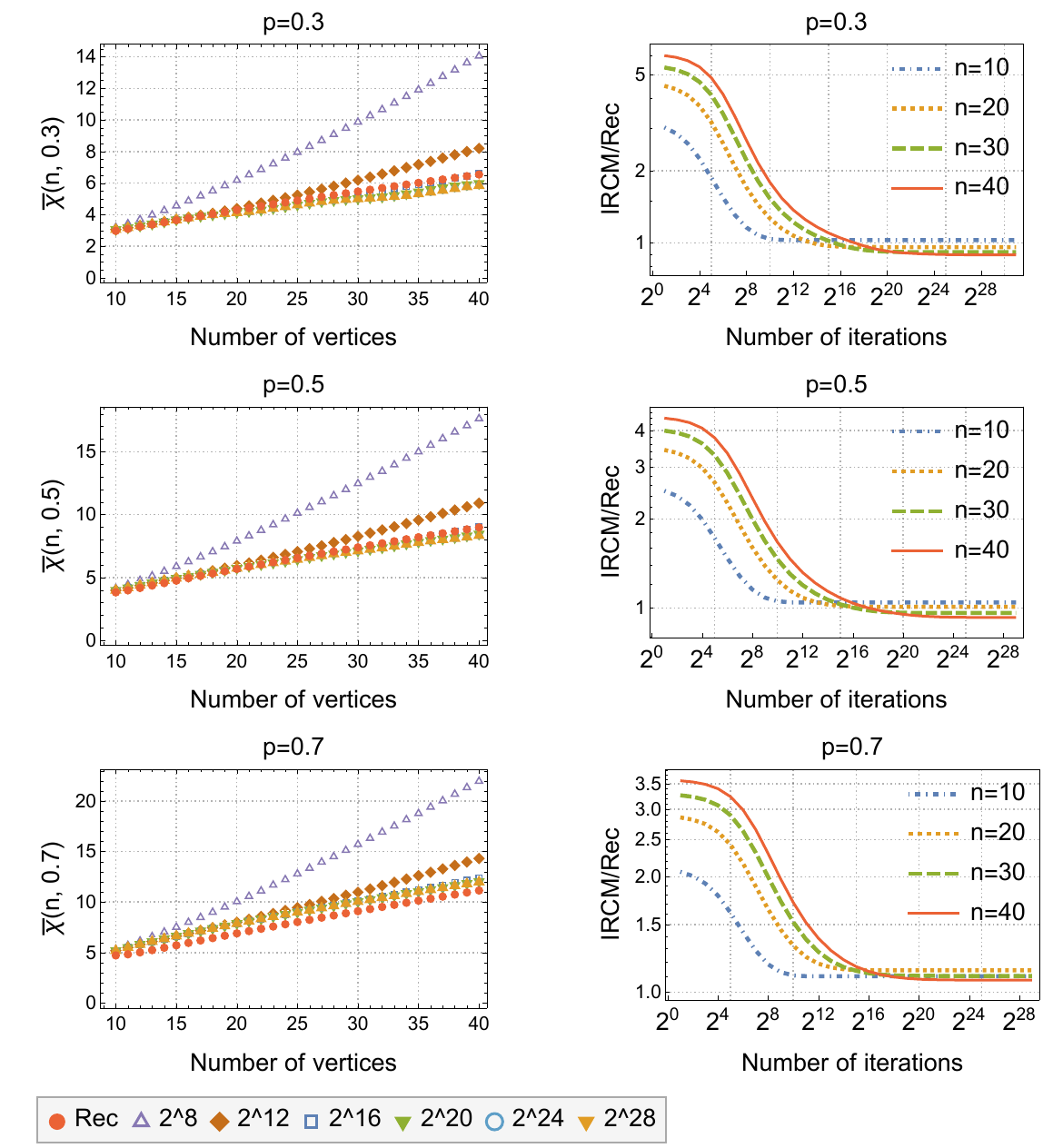} 
   \caption{Comparison of the asymptotic situation of \mcircm \ and the values obtained using the recurrence relation.
Left: Expected values of chromatic number obtained using the recurrence relation (the filled red circles) and expected values of numbers of colors obtained by \mcircm \ when the numbers of iterations are $2^8,2^{12},2^{16},2^{20},2^{24},$ and $2^{28}$.
Right: Ratio of the value obtained by \mcircm \ to that obtained using the recurrence relation (IRCM/Rec). 
Top: $p=0.3$, Middle: $p=0.5$, Bottom: $p=0.7$.
} 
\label{fig:fig_conv1}
\end{figure}

\Cref{fig:fig_conv1} shows a comparison of the asymptotic situation of \mcircm \ and the values obtained using the recurrence relation for $p=0.3, 0.5$, and $0.7$, where the number of samples in the \mc \ method was 16384. 
The left figures of \Cref{fig:fig_conv1} show the asymptotic approach of the values examined by \mcircm \ with increasing number of iterations; the values approach those obtained using the recurrence relation. 
For $p=0.3$ and $0.5$, as the number of \mcircm \ iterations is increased, the \mcircm \ results approach the results obtained using the recurrence relation up to about $2^{16}$ iterations. However, with further iterations, the larger the number of vertices, the smaller the \mcircm \ results than the results obtained using the recurrence relation. For $p=0.7$, the results obtained using the recurrence relation are smaller than the \mcircm \ results. 
Nevertheless, the difference is so small that the values obtained using the recurrence relation are expected to be close to those obtained by \mcircm. 
The right figures of \Cref{fig:fig_conv1} show the asymptotic approach of the ratio of the value obtained by \mcircm \ to that obtained using the recurrence relation, and the final ratio is within the range of $1.0 \pm 0.1$. Thus, the value obtained using the recurrence relation is also expected to be close to that obtained by \mcircm. 

\subsubsection{In the case of graphs with up to $50$ vertices}
\begin{figure}[tb]
 \centering
 \includegraphics[width=\linewidth]{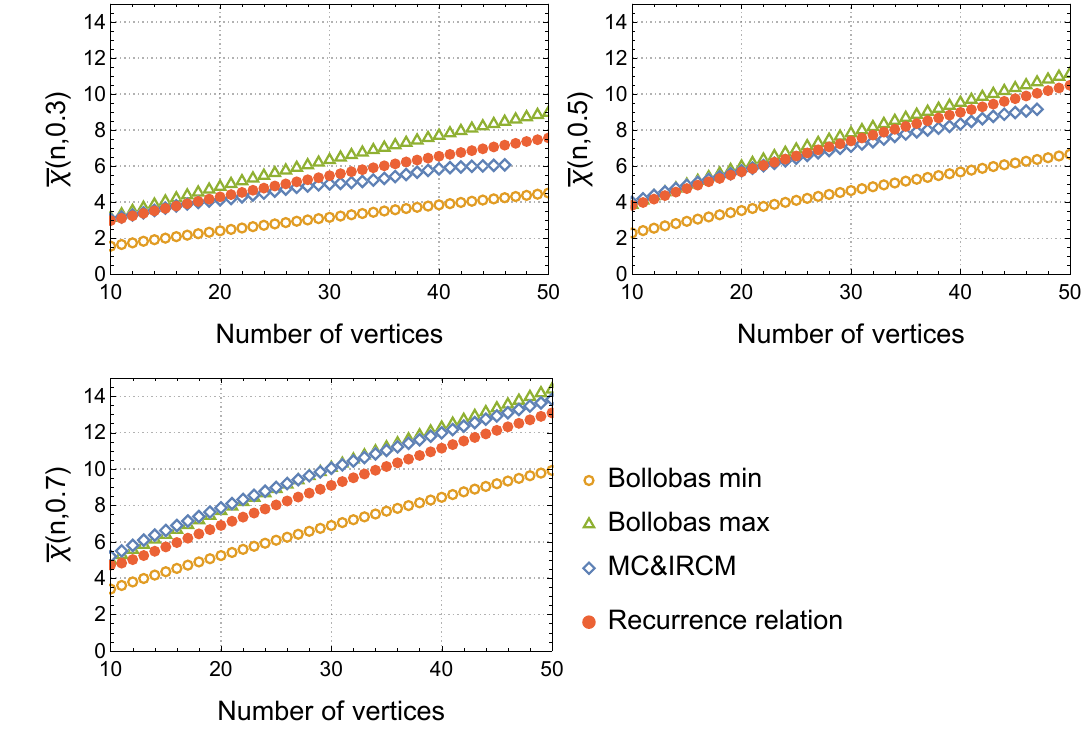} 
 \caption{$\kitaichi{n}{p}$ obtained using the recurrence relation ($n \le 50$). 
 Top left: $p=0.3$, Top right: $p=0.5$, Bottom left: $p=0.7$.
}
 \label{fig:fig_04}
\end{figure}

For relatively small graphs with up to $50$ vertices, the results obtained using the recurrence relation are shown in
\Cref{fig:fig_04}, 
where $p=0.3,0.5,$ and $0.7$, and the initial term of the relation is the exact value at order 9. 
In this figure, the values obtained using the recurrence relation are compared with the experimentally predicted values obtained by \mcircm \ and the values obtained by Bollob\'{a}s.
Note that the values obtained by \agac \ are exact values, and the experimentally predicted values obtained by \mcac \ and \mcircm \ are considered very close to the exact values. Thus, the values obtained by \mcircm \ are also referred to here as exact values.
  
First, we compare the results obtained using the recurrence relation with the values obtained by \mcircm. From the top figures of \Cref{fig:fig_04}, for $p=0.3$ and $0.5$, several values from the initial term are very close to the exact values, but the deviations increase when the order exceeds about $20$ for $p=0.3$ and when the order exceeds about $30$ for $p=0.5$. The bottom left figure of \Cref{fig:fig_04} shows that for $p=0.7$, the values obtained using the recurrence relation already deviate from the exact values at order 10, and the deviation remains almost the same at order $40$. For graphs with an order greater than $40$, the values obtained using the recurrence relation tend to be slightly closer to the exact values.

Next, we compare the results obtained using the recurrence relation with the values obtained using Bollob\'{a}s' relation. From the top figures of \Cref{fig:fig_04}, for $p=0.3$ and $0.5$, the results obtained using the recurrence relation are between the upper bounds of Bollob\'{a}s and the exact values. Thus, in the small graphs where the order is less than $50$, the results obtained using the recurrence relation are closer to the exact values than the values obtained by Bollob\'{a}s. This is because Bollob\'{a}s' estimate is based on the assumption that the number of vertices is very large. From the bottom left figure of \Cref{fig:fig_04}, for $p=0.7$, the values of the upper bounds of Bollob\'{a}s are very close to the exact values. Thus, the values obtained using the recurrence relation are worse estimates than the upper bounds of Bollob\'{a}s. However, taking into account that the values obtained by \mcircm \ are slightly larger than the exact values, the values obtained using the recurrence relation may be closer to the exact values than the upper bounds of Bollob\'{a}s.

\subsubsection{In the case of graphs with very large number of vertices}
\begin{figure}[tb]
 \centering
 \includegraphics[width=\linewidth]{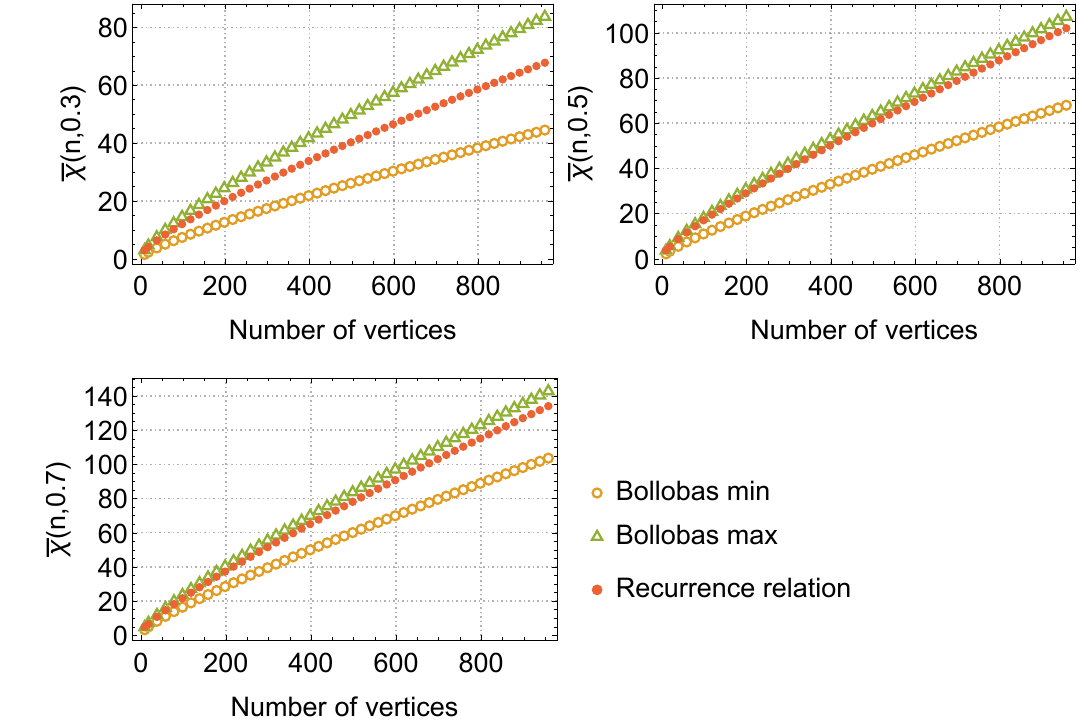} 
 \caption{$\kitaichi{n}{p}$ obtained using the recurrence relation ($n > 50$).  Top left: $p=0.3$, Top right: $p=0.5$, Bottom left: $p=0.7$.}
 \label{fig:fig_05}
\end{figure}

We compared the results obtained using the recurrence relation with Bollob\'{a}s' relation \eqref{eq:bollobas1} for graphs with a large number of vertices, when $p=0.3, 0.5$, and $0.7$. The results are shown in 
\Cref{fig:fig_05}.

This figure shows that the expected value of the chromatic number obtained using the stochastic recurrence relation is slightly smaller than the upper bound of Bollob\'{a}s' relation \eqref{eq:bollobas1}.
However, Bollob\'{a}s' work shows that the chromatic number of almost every random graph is equal to the lower bound in relation \eqref{eq:bollobas1} at the number of vertices $n \to \infty$.
Therefore, as the number of vertices increases, the expected values of the exact chromatic numbers are expected to approach the lower bounds. However, \Cref{fig:fig_05} shows no tendency for the values obtained using the stochastic recurrence relation to asymptotically approach the lower bounds in Bollob\'{a}s' relation (1). Thus, the estimate obtained using the recurrence relation for large-order graphs may be worse than Bollob\'{a}s' lower bound.

\section{Conclusions}
\label{sec:Con}
We attempted to find the expected values of the chromatic numbers of binomial random graphs. For those of small-order graphs, the exact values were obtained by examining all color patterns for each of all the graph patterns (\agac). For large-order graphs, it was difficult to examine all the graph patterns in realistic time. Therefore, the Monte Carlo (\mc) method was used to investigate chromatic numbers to reduce the number of graph patterns to be examined. Furthermore, the larger the number of vertices, the harder it became to examine the entire chromatic pattern in realistic time. Thus, the iterated random color matching method (\mcircm) was used to search chromatic numbers. The experimental predictions of the expected values of the chromatic number obtained in these searches were considered close to the exact values. 

However, for even larger-order random graphs, it became extremely difficult to obtain the expected value of the chromatic number in realistic time. We thus proposed a stochastic recurrence relation method to find the expected value of the chromatic number. 
By comparing the values obtained using the recurrence relation with exact values and experimental predictions, we found that for the probability $p=0.3$ or $0.5$, the results for small graphs with up to about $30$ vertices were closer to the exact values than those of Bollob\'{a}s' previous study. For $p=0.7$, the upper bounds of Bollob\'{a}s' relation were closer to the exact value than the values obtained using the recurrence relation for $n<50$. The possibility that the recurrence relation may approach the exact value in the range $n>=50$ cannot be ruled out, but this could not be clarified in this work.

For very large graphs with more than 100 vertices, the expected values of the chromatic number computed using the recurrence relation were closer to the upper bounds than to the lower bounds of Bollob\'{a}s' relation. Considering that Bollob\'{a}s proved that for $n \to \infty$, almost all random graphs approach Bollob\'{a}s' lower bound, the results obtained using the recurrence relation may be worse than Bollob\'{a}s' estimate. 

To estimate the expected value of the chromatic number with the recurrence relation, we assumed the number of edges connected to the newly added vertex and the number of colors used at the vertices connected by the new edges. The recurrence relation could be improved to obtain values closer to the exact values by finding ways to improve the accuracy of this approximation. 

%% Acknowledgements
\section*{Acknowledgements}
We are grateful to Tsuyoshi Miezaki for his helpful comments and suggestions.

\section*{Funding}
This research did not receive any specific grant from funding agencies in the public, commercial, or not-for-profit sectors.

\appendix
%% This setting is necessary for \Cref.
\renewcommand{\thesection}{\Alph{section}}
\section{Results of all graphs and all coloring search(\agac)}
\label{sec:A}
For random graphs with $2$ to $9$ vertices, %Tables~\ref{tb:exact_edge2}--\ref{tb:exact_edge9} 
\Cref{tb:exact_edge2,tb:exact_edge3,tb:exact_edge4,tb:exact_edge5,tb:exact_edge6,tb:exact_edge7,tb:exact_edge8,tb:exact_edge9} 
show the results of examining all color patterns in all graph patterns to find the chromatic number. These tables show the number of graphs with each chromatic number. These tables also include a breakdown by the number of edges. In these tables, $n$ is the order and $|E|$ is the number of edges of the graph.

\begin{table}[H]  
  \begin{center}
  \caption{Exact number of graphs with each chromatic number ($n=2$)}
  \label{tb:exact_edge2}
  %\begin{threeparttable}[h]
    \setlength{\tabcolsep}{10pt}
    \begin{tabular}{ccccc} 
    \hline
    & \multicolumn{2}{c}{Chromatic number}  \\
    \cline{2-3} 
    %$|E|$ \tnote{$\dagger$}  &1 & 2& Total\\
    $|E|$ &1 & 2& Total\\
    \hline
    0 & 1 & 0  & 1 \\ 
    1 & 0 & 1  & 1 \\ 
    \hline
    Total & 1 & 1 & 2 \\
    \hline 
    \end{tabular}
%    \begin{tablenotes}
%      \item[$\dagger$] Number of edges.
%    \end{tablenotes}
%    \end{threeparttable}
  \end{center}
\end{table}%

\begin{table}[H]
  \begin{center}
  \caption{Exact number of graphs with each chromatic number ($n=3$)}
  \label{tb:exact_edge3}
%  \begin{threeparttable}[h]
    \setlength{\tabcolsep}{10pt}
    \begin{tabular}{ccccc} 
    \hline
    & \multicolumn{3}{c}{Chromatic number} \\
    \cline{2-4} 
%    $|E|$ \tnote{$\dagger$}  & 1 & 2 & 3 & Total\\
    $|E|$ & 1 & 2 & 3 & Total\\
    \hline
    0 & 1 & 0 & 0  & 1 \\ 
    1 & 0 & 3 & 0  & 3 \\ 
    2 & 0 & 3 & 0  & 3  \\ 
    3 & 0 & 0 & 1  & 1 \\
    \hline
    Total & 1 & 6 & 1 & 8\\
    \hline 
    \end{tabular}
%    \begin{tablenotes}
%      \item[$\dagger$] Number of edges.
%    \end{tablenotes}
%    \end{threeparttable}
  \end{center}
\end{table}%

\begin{table}[H]
  \begin{center}
  \caption{Exact number of graphs with each chromatic number ($n=4$)}
  \label{tb:exact_edge4}
%  \begin{threeparttable}[h]
  \setlength{\tabcolsep}{10pt}
    \begin{tabular}{cccccc} 
    \hline
    & \multicolumn{4}{c}{Chromatic number}  \\
    \cline{2-5} 
%    $|E|$ \tnote{$\dagger$}  &1 & 2& 3 & 4 & Total\\
    $|E|$ &1 & 2& 3 & 4 & Total\\
    \hline
    0 & 1 & 0 & 0 & 0 & 1 \\ 
    1 & 0 & 6 & 0 & 0 & 6  \\ 
    2 & 0 & 15 & 0 & 0 & 15  \\ 
    3 & 0 & 16 & 4 & 0 & 20  \\ 
    4 & 0 & 3 & 12 & 0 & 15  \\ 
    5 & 0 & 0 & 6 & 0 & 6  \\ 
    6 & 0 & 0 & 0 & 1 & 1  \\
    \hline 
    Total & 1 & 40 & 22 & 1  & 64 \\
    \hline 
    \end{tabular}
%    \begin{tablenotes}
%      \item[$\dagger$] Number of edges.
%    \end{tablenotes}
%    \end{threeparttable}
  \end{center}
\end{table}%

\begin{table}[H]
  \begin{center}
  \caption{Exact number of graphs with each chromatic number ($n=5$)}
  \label{tb:exact_edge5}
%  \begin{threeparttable}[h]
  \setlength{\tabcolsep}{10pt}
    \begin{tabular}{ccccccc} 
    \hline
    & \multicolumn{5}{c}{Chromatic number}  \\
    \cline{2-6} 
%    $|E|$ \tnote{$\dagger$}  &1 & 2& 3 & 4 & 5 & Total\\
    $|E|$ &1 & 2& 3 & 4 & 5 & Total\\
    \hline
    0 & 1 & 0 & 0 & 0 & 0 & 1 \\
    1 & 0 & 1 & 0 & 0 & 0 & 10 \\
    2 & 0 & 45 & 0 & 0 & 0 & 45 \\
    3 & 0 & 110 & 10 & 0 & 0 & 120 \\
    4 & 0 & 140 & 70 & 0 & 0 & 210 \\
    5 & 0 & 60 & 	192 & 0 & 0 & 252 \\
    6 & 0 & 10 & 195 & 5 & 0 & 210 \\
    7 & 0 & 0 & 100 & 20 & 0 & 120 \\
    8 & 0 & 0 & 15 & 30 & 0 & 45 \\
    9 & 0 & 0 & 0 & 10 & 0 & 10 \\
    10 & 0 & 0 & 0 & 0 & 1 & 1 \\
    \hline 
    Total & 1 & 375 & 582 & 65 & 1 & 1024 \\
    \hline 
    \end{tabular}
%    \begin{tablenotes}
%      \item[$\dagger$] Number of edges.
%    \end{tablenotes}
%    \end{threeparttable}
  \end{center}
\end{table}%

\begin{table}[H]
  \begin{center}
  \caption{Exact number of graphs with each chromatic number ($n=6$)}
  \label{tb:exact_edge6}
  %\begin{threeparttable}[h]
  \setlength{\tabcolsep}{10pt}
    \begin{tabular}{cccccccc} 
    \hline
    & \multicolumn{6}{c}{Chromatic number}  \\
    \cline{2-7} 
%    $|E|$ \tnote{$\dagger$}  &1 & 2& 3 & 4 & 5 & 6 & Total\\
    $|E|$ &1 & 2& 3 & 4 & 5 & 6 & Total\\
    \hline
    0 & 1 & 0 & 0 & 0 & 0 & 0 & 1 \\
    1 & 0 & 15 & 0 & 0 & 0 & 0 & 15 \\
    2 & 0 & 105 & 0 & 0 & 0 & 0 & 105 \\
    3 & 0 & 435 & 20 & 0 & 0 & 0 & 455 \\
    4 & 0 & 1125 & 240 & 0 & 0 & 0 & 1365 \\
    5 & 0 & 1701 & 1302 & 0 & 0 & 0 & 3003 \\
    6 & 0 & 1200 & 3790 & 15 & 0 & 0 & 5005 \\
    7 & 0 & 480 & 5820 & 135 & 0 & 0 & 6435 \\
    8 & 0 & 105 & 5790 & 540 & 0 & 0 & 6435 \\
    9 & 0 & 10 & 3795 & 1200 & 0 & 0 & 5005 \\
   10 & 0 & 0 & 1365 & 1632 & 6 & 0 & 3003 \\
   11 & 0 & 0 & 240 & 1095 & 30 & 0 & 1365 \\
   12 & 0 & 0 & 15 & 380 & 60 & 0 & 455 \\
   13 & 0 & 0 & 0 & 45 & 60 & 0 & 105 \\
   14 & 0 & 0 & 0 & 0 & 15 & 0 & 15 \\
   15 & 0 & 0 & 0 & 0 & 0 & 1 & 1 \\
    \hline 
    Total & 1 & 5176 & 22377 & 5042 & 171 & 1 & 32768 \\
    \hline 
    \end{tabular}
%    \begin{tablenotes}
%      \item[$\dagger$] Number of edges.
%    \end{tablenotes}
%    \end{threeparttable}
  \end{center}
\end{table}%

\begin{table}[H]  
  \begin{center}
  \caption{Exact number of graphs with each chromatic number ($n=7$)}
  \label{tb:exact_edge7}
%  \begin{threeparttable}[h]
  \setlength{\tabcolsep}{10pt}
    \begin{tabular}{ccccccccc} 
    \hline
    & \multicolumn{7}{c}{Chromatic number}  \\
    \cline{2-8} 
%    $|E|$ \tnote{$\dagger$}  &1 & 2& 3 & 4 & 5 & 6 & 7 & Total\\
    $|E|$ &1 & 2& 3 & 4 & 5 & 6 & 7 & Total\\
    \hline
0 & 1 & 0 & 0 & 0 & 0 & 0 & 0 & 1 \\
1 & 0 & 21 & 0 & 0 & 0 & 0 & 0 & 21 \\
2 & 0 & 210 & 0 & 0 & 0 & 0 & 0 & 210 \\
3 & 0 & 1295 & 35 & 0 & 0 & 0 & 0 & 1330 \\
4 & 0 & 5355 & 630 & 0 & 0 & 0 & 0 & 5985 \\
5 & 0 & 14952 & 5397 & 0 & 0 & 0 & 0 & 20349 \\
6 & 0 & 26572 & 27657 & 35 & 0 & 0 & 0 & 54264 \\
7 & 0 & 26670 & 89085 & 525 & 0 & 0 & 0 & 116280 \\
8 & 0 & 17535 & 182280 & 3675 & 0 & 0 & 0 & 203490 \\
9 & 0 & 7840 & 270375 & 15715 & 0 & 0 & 0 & 293930 \\
10 & 0 & 2331 & 304500 & 45864 & 21 & 0 & 0 & 352716 \\
11 & 0 & 420 & 256326 & 95739 & 231 & 0 & 0 & 352716 \\
12 & 0 & 35 & 147630 & 145110 & 1155 & 0 & 0 & 293930 \\
13 & 0 & 0 & 55125 & 144900 & 3465 & 0 & 0 & 203490 \\
14 & 0 & 0 & 12915 & 96540 & 6825 & 0 & 0 & 116280 \\
15 & 0 & 0 & 1750 & 43477 & 9030 & 7 & 0 & 54264 \\
16 & 0 & 0 & 105 & 12075 & 8127 & 42 & 0 & 20349 \\
17 & 0 & 0 & 0 & 1785 & 4095 & 105 & 0 & 5985 \\
18 & 0 & 0 & 0 & 105 & 1085 & 140 & 0 & 1330 \\
19 & 0 & 0 & 0 & 0 & 105 & 105 & 0 & 210 \\
20 & 0 & 0 & 0 & 0 & 0 & 21 & 0 & 21 \\
21 & 0 & 0 & 0 & 0 & 0 & 0 & 1 & 1 \\
    \hline 
    Total & 1 & 103236 & 1353810 & 605545 & 34139 & 420 & 1 & 2097152 \\
    \hline 
    \end{tabular}
%    \begin{tablenotes}
%      \item[$\dagger$] Number of edges.
%    \end{tablenotes}
%    \end{threeparttable}
  \end{center}
\end{table}%

\begin{table}[H]
  \begin{center}
  \caption{Exact number of graphs with each chromatic number ($n=8$)}
  \label{tb:exact_edge8}
%  \begin{threeparttable}[h]
  \setlength{\tabcolsep}{3pt}
    \begin{tabular}{cccccccccc} 
    \hline
    & \multicolumn{8}{c}{Chromatic number}  \\
    \cline{2-9} 
%    $|E|$ \tnote{$\dagger$} &1 & 2& 3 & 4 & 5 & 6 & 7 & 8 & Total\\
    $|E|$ &1 & 2& 3 & 4 & 5 & 6 & 7 & 8 & Total\\
    \hline
0 & 1 & 0 & 0 & 0 & 0 & 0 & 0 & 0 & 1 \\
1 & 0 & 28 & 0 & 0 & 0 & 0 & 0 & 0 & 28 \\
2 & 0 & 378 & 0 & 0 & 0 & 0 & 0 & 0 & 378 \\
3 & 0 & 3220 & 56 & 0 & 0 & 0 & 0 & 0 & 3276 \\
4 & 0 & 19075 & 1400 & 0 & 0 & 0 & 0 & 0 & 20475 \\
5 & 0 & 81228 & 17052 & 0 & 0 & 0 & 0 & 0 & 98280 \\
6 & 0 & 246414 & 130256 & 70 & 0 & 0 & 0 & 0 & 376740 \\
7 & 0 & 507424 & 675076 & 1540 & 0 & 0 & 0 & 0 & 1184040 \\
8 & 0 & 666015 & 2425920 & 16170 & 0 & 0 & 0 & 0 & 3108105 \\
9 & 0 & 620900 & 6178760 & 107240 & 0 & 0 & 0 & 0 & 6906900 \\
10 & 0 & 431368 & 12187980 & 503706 & 56 & 0 & 0 & 0 & 13123110 \\
11 & 0 & 226296 & 19464228 & 1782648 & 1008 & 0 & 0 & 0 & 21474180 \\
12 & 0 & 88928 & 25388076 & 4936183 & 8568 & 0 & 0 & 0 & 30421755 \\
13 & 0 & 25480 & 26528068 & 10842916 & 45696 & 0 & 0 & 0 & 37442160 \\
14 & 0 & 5040 & 21149520 & 18791100 & 170940 & 0 & 0 & 0 & 40116600 \\
15 & 0 & 616 & 12471620 & 24495716 & 474180 & 28 & 0 & 0 & 37442160 \\
16 & 0 & 35 & 5413828 & 24000872 & 1006656 & 364 & 0 & 0 & 30421755 \\
17 & 0 & 0 & 1725528 & 18084276 & 1662192 & 2184 & 0 & 0 & 21474180 \\
18 & 0 & 0 & 396130 & 10575096 & 2143876 & 8008 & 0 & 0 & 13123110 \\
19 & 0 & 0 & 62440 & 4669140 & 2155300 & 20020 & 0 & 0 & 6906900 \\
20 & 0 & 0 & 6090 & 1480227 & 1585920 & 35868 & 0 & 0 & 3108105 \\
21 & 0 & 0 & 280 & 319200 & 817736 & 46816 & 8 & 0 & 1184040 \\
22 & 0 & 0 & 0 & 43750 & 288806 & 44128 & 56 & 0 & 376740 \\
23 & 0 & 0 & 0 & 3360 & 65100 & 29652 & 168 & 0 & 98280 \\
24 & 0 & 0 & 0 & 105 & 8120 & 11970 & 280 & 0 & 20475 \\
25 & 0 & 0 & 0 & 0 & 420 & 2576 & 280 & 0 & 3276 \\
26 & 0 & 0 & 0 & 0 & 0 & 210 & 168 & 0 & 378 \\
27 & 0 & 0 & 0 & 0 & 0 & 0 & 28 & 0 & 28 \\
28 & 0 & 0 & 0 & 0 & 0 & 0 & 0 & 1 & 1 \\
    \hline 
    Total & 1 & 2922445 & 134222308 & 120653315 & 10434574 & 201824 & 988 & 1 & 268435456 \\
    \hline 
    \end{tabular}
%    \begin{tablenotes}
%      \item[$\dagger$] Number of edges.
%    \end{tablenotes}
%    \end{threeparttable}
  \end{center}
\end{table}%

\begin{table}[H]
  \begin{center}
  \caption{Exact number of graphs with each chromatic number ($n=9$)}
  \label{tb:exact_edge9}
%  \begin{threeparttable}[h]
  \setlength{\tabcolsep}{3pt}
  \scriptsize
    \begin{tabular}{ccccccccccc} 
    \hline
    & \multicolumn{9}{c}{Chromatic number}  \\
    \cline{2-10} 
%    $|E|$ \tnote{$\dagger$}  &1 & 2& 3 & 4 & 5 & 6 & 7 & 8 & 9 & Total\\
    $|E|$ &1 & 2& 3 & 4 & 5 & 6 & 7 & 8 & 9 & Total\\
    \hline
0 & 1 & 0 & 0 & 0 & 0 & 0 & 0 & 0 & 0 & 1 \\
1 & 0 & 36 & 0 & 0 & 0 & 0 & 0 & 0 & 0 & 36 \\
2 & 0 & 630 & 0 & 0 & 0 & 0 & 0 & 0 & 0 & 630 \\
3 & 0 & 7056 & 84 & 0 & 0 & 0 & 0 & 0 & 0 & 7140 \\
4 & 0 & 56133 & 2772 & 0 & 0 & 0 & 0 & 0 & 0 & 58905 \\
5 & 0 & 331884 & 45108 & 0 & 0 & 0 & 0 & 0 & 0 & 376992 \\
6 & 0 & 1475964 & 471702 & 126 & 0 & 0 & 0 & 0 & 0 & 1947792 \\
7 & 0 & 4864680 & 3479220 & 3780 & 0 & 0 & 0 & 0 & 0 & 8347680 \\
8 & 0 & 11445534 & 18759996 & 54810 & 0 & 0 & 0 & 0 & 0 & 30260340 \\
9 & 0 & 18626412 & 75006568 & 510300 & 0 & 0 & 0 & 0 & 0 & 94143280 \\
10 & 0 & 22709736 & 228051306 & 3425688 & 126 & 0 & 0 & 0 & 0 & 254186856 \\
11 & 0 & 21711060 & 561415176 & 17675784 & 3276 & 0 & 0 & 0 & 0 & 600805296 \\
12 & 0 & 16618140 & 1161987645 & 73030965 & 40950 & 0 & 0 & 0 & 0 & 1251677700 \\
13 & 0 & 10257156 & 2052036504 & 248168340 & 327600 & 0 & 0 & 0 & 0 & 2310789600 \\
14 & 0 & 5095476 & 3084952374 & 704366910 & 1882440 & 0 & 0 & 0 & 0 & 3796297200 \\
15 & 0 & 2016504 & 3883815216 & 1673809440 & 8261316 & 84 & 0 & 0 & 0 & 5567902560 \\
16 & 0 & 623007 & 3976618527 & 3301885944 & 28742868 & 1764 & 0 & 0 & 0 & 7307872110 \\
17 & 0 & 145152 & 3244054968 & 5272045380 & 81233460 & 17640 & 0 & 0 & 0 & 8597496600 \\
18 & 0 & 24024 & 2102745582 & 6782676264 & 189577710 & 111720 & 0 & 0 & 0 & 9075135300 \\
19 & 0 & 2520 & 1088934588 & 7138556712 & 369500040 & 502740 & 0 & 0 & 0 & 8597496600 \\
20 & 0 & 126 & 452161962 & 6248235987 & 605765475 & 1708560 & 0 & 0 & 0 & 7307872110 \\
21 & 0 & 0 & 150141096 & 4575928248 & 837286596 & 4546584 & 36 & 0 & 0 & 5567902560 \\
22 & 0 & 0 & 39441780 & 2773801674 & 973368846 & 9684360 & 540 & 0 & 0 & 3796297200 \\
23 & 0 & 0 & 8028216 & 1357270992 & 928757844 & 16728768 & 3780 & 0 & 0 & 2310789600 \\
24 & 0 & 0 & 1222263 & 521127327 & 705721716 & 23590014 & 16380 & 0 & 0 & 1251677700 \\
25 & 0 & 0 & 131040 & 153308988 & 420147504 & 27168624 & 49140 & 0 & 0 & 600805296 \\
26 & 0 & 0 & 8820 & 33801390 & 194865048 & 25403490 & 108108 & 0 & 0 & 254186856 \\
27 & 0 & 0 & 280 & 5414220 & 69458760 & 19090092 & 179928 & 0 & 0 & 94143280 \\
28 & 0 & 0 & 0 & 595980 & 18400707 & 11034072 & 229572 & 9 & 0 & 30260340 \\
29 & 0 & 0 & 0 & 40320 & 3439296 & 4643964 & 224028 & 72 & 0 & 8347680 \\
30 & 0 & 0 & 0 & 1260 & 422310 & 1359918 & 164052 & 252 & 0 & 1947792 \\
31 & 0 & 0 & 0 & 0 & 30240 & 258300 & 87948 & 504 & 0 & 376992 \\
32 & 0 & 0 & 0 & 0 & 945 & 27720 & 29610 & 630 & 0 & 58905 \\
33 & 0 & 0 & 0 & 0 & 0 & 1260 & 5376 & 504 & 0 & 7140 \\
34 & 0 & 0 & 0 & 0 & 0 & 0 & 378 & 252 & 0 & 630 \\
35 & 0 & 0 & 0 & 0 & 0 & 0 & 0 & 36 & 0 & 36 \\
36 & 0 & 0 & 0 & 0 & 0 & 0 & 0 & 0 & 1 & 1 \\
    \hline 
    Total & 1 & 116011230 & 22133512793 & 40885736829 & 5437235073 & 145879674 & 1098876 & 2259 & 1 & 68719476736 \\
    \hline
    \end{tabular}
%    \begin{tablenotes}
%      \item[$\dagger$] Number of edges.
%    \end{tablenotes}
%    \end{threeparttable} 
  \end{center}  
\end{table}%

\section{Results of Monte Carlo method and all coloring search (MC\&AC)}
\label{sec:B}
\setcounter{table}{0}
\Cref{tb:monte_all} shows the results obtained by \mcac \ for the $16384$ graphs extracted by the \mc \ method with probabilities $p = 0.3, 0.5$, and $0.7$. Only results for which all chromatic patterns could be investigated in realistic time are shown. 
\begin{table}[H]
  \caption{Expected values of the chromatic number obtained by \mcac \ (16384 samples)}
  \label{tb:monte_all}
  \begin{center}
  \setlength{\tabcolsep}{15pt}
  {
    \begin{tabular}{cccc} 
    \hline
    $n$  & $p=0.3$ & $p=0.5$ & $p=0.7$ \\
    \hline
    10 & 3.05609 & 3.98541 & 5.18646  \\
    11 & 3.17175 & 4.18091 & 5.50482  \\
    12 & 3.28790 & 4.37830 & 5.78979  \\
    13 & 3.41375 & 4.56677 & 6.08941  \\
    14 & 3.53948 & 4.76459 &   \\
    15 & 3.67120 &  &   \\
    \hline      
    \end{tabular}
  }
  \end{center}
\end{table}%

\section{Results of Monte Carlo method and IRCM method (\mcircm)}
\label{sec:C}
\setcounter{table}{0}
\Cref{tb:monte_IRCM} shows the results obtained by the \ircm \ method for the $16384$ graphs extracted by the \mc \ method with probabilities $p = 0.3, 0.5$, and $0.7$. Only results for graphs where the predicted value of the chromatic number was obtained in realistic time are shown. 
\begin{table}[H]
  \caption{Expected values of chromatic number obtained by \mcircm \ (16384 samples)}
  \label{tb:monte_IRCM}
  \begin{center}
  \setlength{\tabcolsep}{10pt}
  {
    \begin{tabular}{ccccccccc} 
    \cline{1-4}
    \cline{6-9}
    $n$  & $p=0.3$ & $p=0.5$ & $p=0.7$ & & $n$  & $p=0.3$ & $p=0.5$ & $p=0.7$\\
    \cline{1-4}
    \cline{6-9}
    7  & 2.57208 & 3.28235  & 4.14459 & & 29 & 4.96112 & 6.99377 & 9.83837 \\
    8  & 2.75994 & 3.51740 & 4.50994 & & 30 & 5.00939 & 7.10309 & 10.04406 \\
    9  & 2.92694 & 3.76068 & 4.85943 & & 31 & 5.04583 & 7.22601 & 10.22991 \\
    10 & 3.06158 & 3.98566 & 5.19287 & & 32 & 5.08734 & 7.34644 & 10.44055 \\
    11 & 3.17449 & 4.18341 & 5.49645 & & 33 & 5.14721 & 7.49286 & 10.63226 \\
    12 & 3.28973 & 4.37543 & 5.80181 & & 34 & 5.22467 & 7.63953 & 10.83416 \\
    13 & 3.41223 & 4.57300 & 6.09307 & & 35 & 5.31707 & 7.78149 & 11.02697 \\
    14 & 3.54333 & 4.76019 & 6.36535 & & 36 & 5.42236 & 7.90857 & 11.22503 \\
    15 & 3.67327 & 4.93408 & 6.62200 & & 37 & 5.53637 & 8.00928 & 11.41369 \\
    16 & 3.79565 & 5.10516 & 6.89300 & & 38 & 5.65606 & 8.10284 & 11.61236 \\
    17 & 3.90118 & 5.25275 & 7.13995 & & 39 & 5.76123 & 8.20667 & 11.80828 \\
    18 & 3.99200 & 5.41376 & 7.38824 & & 40 & 5.84442 & 8.33148 &11.99114 \\
    19 & 4.06359 & 5.57428 & 7.63586 & & 41 & 5.91729 & 8.47058 & 12.17205 \\
    20 & 4.13464 & 5.73358 & 7.86962 & & 42 & 5.95379 & 8.62299 & 12.36236 \\
    21 & 4.20379 & 5.89117 & 8.09716 & & 43 & 5.98645 & 8.76917 & 12.55267 \\
    22 & 4.29089 & 6.03149 & 8.32849 & & 44 & 6.00903 & 8.88422 & 12.74066 \\
    23 & 4.39971 & 6.16010 & 8.55572 & & 45 & 6.03997 & 8.97583 & 12.92150 \\
    24 & 4.50317 & 6.28992 & 8.78234 & & 46 & 6.06848 & 9.06122 & 13.09875 \\
    25 & 4.61462 & 6.43842 & 8.98883 & & 47 &  & 9.14917 & 13.27313 \\
    26 & 4.72460 & 6.58984 & 9.20043 & & 48 &  & 9.25470 & 13.46356 \\
    27 & 4.82781 & 6.73486 & 9.42053 & & 49 &  & 9.38312 & 13.64019 \\
    28 & 4.89990 & 6.87140 & 9.63208 & & 50 &  &  & 13.83203 \\
    \cline{1-4}
    \cline{6-9}     
    \end{tabular}
  }
  \end{center}
\end{table}%

%%       \cite{<label>} ==> [1]

%%
%Example citation, See \cite{lamport94}.

\bibliographystyle{elsarticle-num} 
\bibliography{citation}




%% For numbered reference style

%\bibitem{lamport94}
%  \textit{\LaTeX: a document preparation system},
%  Addison Wesley, Massachusetts,

%\end{thebibliography}
\end{document}